\documentclass{article}

\usepackage{amssymb}

\newcommand{\newc}{\newcommand}

\newc{\eqnoset}{\setcounter{equation}{0}}

\newcommand{\mref}[1]{(\ref{#1})}

\newcommand{\reftheo}[1]{Theorem~\ref{#1}}

\newcommand{\refsec}[1]{Section~\ref{#1}}

\newcommand{\beq}{\begin{equation}}
	\newcommand{\eeq}{\end{equation}}
\newcommand{\beqno}[1]{\begin{equation}\label{#1}}
	
	\newcommand{\barr}{\begin{array}}
		\newcommand{\earr}{\end{array}}
	
	\newc{\bearr}{\begin{eqnarray*}}
		\newc{\eearr}{\end{eqnarray*}}
	
	\newc{\bearrno}[1]{\begin{eqnarray}\label{#1}}
		\newc{\eearrno}{\end{eqnarray}}
	
	\newc{\non}{\nonumber}
	\newc{\nol}{\nonumber\nl}
	
	\newcommand{\bdes}{\begin{description}}
		\newcommand{\edes}{\end{description}}
	\newc{\benu}{\begin{enumerate}}
		\newc{\eenu}{\end{enumerate}}
	\newc{\btab}{\begin{tabular}}
		\newc{\etab}{\end{tabular}}

	\newtheorem{theorem}{Theorem}[section]
	\newtheorem{defi}[theorem]{Definition}
	\newtheorem{lemma}[theorem]{Lemma}
	\newtheorem{rem}[theorem]{Remark}
	\newtheorem{exam}[theorem]{Example}
	\newtheorem{propo}[theorem]{Proposition}
	\newtheorem{corol}[theorem]{Corollary}
	\newtheorem{conj}[theorem]{Conjecture}

	\newcommand{\btheo}[1]{\begin{theorem}\label{#1}}
		\newc{\brem}[1]{\begin{rem}\label{#1}\em}
			\newc{\bexam}[1]{\begin{exam}\label{#1}\em}
				\newc{\bdefi}[1]{\begin{defi}\label{#1}}
					\newcommand{\blemm}[1]{\begin{lemma}\label{#1}}
						\newcommand{\bprop}[1]{\begin{propo}\label{#1}}
							\newcommand{\bcoro}[1]{\begin{corol}\label{#1}}
								\newcommand{\btheoc}[1]{\begin{conj}\label{#1}}
									\newcommand{\etheo}{\end{theorem}}
								\newc{\etheoc}{\end{conj}}
							\newcommand{\elemm}{\end{lemma}}
						\newcommand{\eprop}{\end{propo}}
					\newcommand{\ecoro}{\end{corol}}
				\newc{\erem}{\end{rem}}
			\newc{\eexam}{\end{exam}}
		\newc{\edefi}{\end{defi}}
	
	\newc{\rmk}[1]{{\bf REMARK #1: }}
	\newc{\DN}[1]{{\bf DEFINITION #1: }}
	
	\newcommand{\bproof}{{\bf Proof:~~}}
	\newc{\eproof}{{\vrule height8pt width5pt depth0pt}\vspace{3mm}}
	
	\newc{\bfrac}[2]{\dspl{\frac{#1}{#2}}}

	\newc{\nid}{\noindent}
	
	\newcommand{\dspl}{\displaystyle}
	\newc{\grad}{\nabla}
	\newc{\Div}{\mbox{div}}
	\newc{\pdt}[1]{\dspl{\frac{\partial{#1}}{\partial t}}}
	\newc{\pdn}[1]{\dspl{\frac{\partial{#1}}{\partial \nu}}}
	\newc{\pdNi}[1]{\dspl{\frac{\partial{#1}}{\partial \mathcal{N}_i}}}
	\newc{\pD}[2]{\dspl{\frac{\partial{#1}}{\partial #2}}}
	\newc{\dt}{\dspl{\frac{d}{dt}}}
	\newc{\bdry}[1]{\mbox{$\partial #1$}}
	\newc{\sgn}{\mbox{sign}}
	
	\newc{\Hess}[1]{\frac{\partial^2 #1}{\pdh z_i \pdh z_j}}
	\newc{\hess}[1]{\partial^2 #1/\pdh z_i \pdh z_j}

	\newc{\ag}{\alpha}
	\newc{\bg}{\beta}
	\newc{\cg}{\gamma}\newc{\Cg}{\Gamma}
	\newc{\dg}{\delta}\newc{\Dg}{\Delta}
	\newc{\eg}{\varepsilon}
	\newc{\zg}{\zeta}
	\newc{\thg}{\theta}
	\newc{\llg}{\lambda}\newc{\LLg}{\Lambda}
	\newc{\kg}{\kappa}
	\newc{\rg}{\rho}
	\newc{\sg}{\sigma}\newc{\Sg}{\Sigma}
	\newc{\tg}{\tau}
	\newc{\fg}{\phi}\newc{\Fg}{\Phi}
	\newc{\vfg}{\varphi}
	\newc{\og}{\omega}\newc{\Og}{\Omega}
	\newc{\pdh}{\partial}
	
	\newc{\ccG}{{\cal G}}

	\newc{\ii}[1]{\int_{#1}}
	\newc{\iidx}[2]{{\dspl\int_{#1}~#2~dx}}
	\newc{\bii}[1]{{\dspl \ii{#1} }}
	\newc{\biii}[2]{{\dspl \iii{#1}{#2} }}
	\newc{\su}[2]{\sum_{#1}^{#2}}
	\newc{\bsu}[2]{{\dspl \su{#1}{#2} }}

	\newc{\biiom}[1]{{\dspl\int_{\bdrom}~ #1 ~d\sg}}
	\newc{\io}[1]{{\dspl\int_{\Og}~ #1 ~dx}}
	\newc{\bio}[1]{{\dspl\int_{\bdrom}~ #1 ~d\sg}}
	\newc{\bsir}{\bsu{i=1}{r}}
	\newc{\bsim}{\bsu{i=1}{m}}
	
	\newc{\iibr}[2]{\iidx{\bprw{#1}}{#2}}
	\newc{\Intbr}[1]{\iibr{R}{#1}}
	\newc{\intbr}[1]{\iibr{\rg}{#1}}
	\newc{\intt}[3]{\int_{#1}^{#2}\int_\Og~#3~dxdt}
	\newc{\itQ}[2]{\dspl{\int\hspace{-2.5mm}\int_{#1}~#2~dz}}
	\newc{\mitQ}[2]{\dspl{\rule[1mm]{4mm}{.3mm}\hspace{-5.3mm}\int\hspace{-2.5mm}\int_{#1}~#2~dz}}
	\newc{\mitQQ}[3]{\dspl{\rule[1mm]{4mm}{.3mm}\hspace{-5.3mm}\int\hspace{-2.5mm}\int_{#1}~#2~#3}}
	
	\newc{\mitx}[2]{\dspl{\rule[1mm]{3mm}{.3mm}\hspace{-4mm}\int_{#1}~#2~dx}}
	\newc{\mitmu}[2]{\dspl{\rule[1mm]{3mm}{.3mm}\hspace{-4mm}\int_{#1}~#2~d\mu}}
	\newc{\iidmu}[2]{\iidx{#1}{#2}}
	
	\newc{\iidm}[3]{{\dspl\int_{#1}~#2~d #3}}
	
	\newc{\itQmu}[2]{\dspl{\int\hspace{-2.5mm}\int_{#1}~#2~d\mu}}
	\newc{\mitQmu}[2]{\dspl{\rule[1mm]{4mm}{.3mm}\hspace{-5.3mm}\int\hspace{-2.5mm}\int_{#1}~#2~d\mu}}

	\newc{\mitQq}[2]{\dspl{\rule[1mm]{4mm}{.3mm}\hspace{-5.3mm}\int\hspace{-2.5mm}\int_{#1}~#2~d\bar{z}}}
	\newc{\itQq}[2]{\dspl{\int\hspace{-2.5mm}\int_{#1}~#2~d\bar{z}}}

	\newc{\pder}[2]{\dspl{\frac{\partial #1}{\partial #2}}}

	\newc{\bdrom}{\bdry{\Og}}

	\newc{\bilhom}{\mbox{Bil}(\mbox{Hom}(\RR^{nm},\RR^{nm}))}
	\newc{\VV}[1]{{V(Q_{#1})}}
	
	\newc{\ccA}{{\mathcal A}}
	\newc{\ccB}{{\mathcal B}}
	\newc{\ccC}{{\mathcal C}}
	\newc{\ccD}{{\mathcal D}}
	\newc{\ccE}{{\mathcal E}}
	\newc{\ccH}{\mathcal{H}}
	\newc{\ccF}{\mathcal{F}}
	\newc{\ccI}{{\mathcal I}}
	\newc{\ccJ}{{\mathcal J}}
	\newc{\ccK}{{\mathcal K}}
	\newc{\ccP}{{\mathcal P}}
	\newc{\ccQ}{{\mathcal Q}}
	\newc{\ccR}{{\mathcal R}}
	\newc{\ccS}{{\mathcal S}}
	\newc{\ccT}{{\mathcal T}}
	\newc{\ccX}{{\mathcal X}}
	\newc{\ccY}{{\mathcal Y}}
	\newc{\ccZ}{{\mathcal Z}}
	
	\newc{\bb}[1]{{\mathbf #1}}
	
	\newc{\myprod}[1]{\langle #1 \rangle}
	\newc{\mypar}[1]{\left( #1 \right)}
	\newc{\BLLg}{\mathbf{\LLg}}
	
	\newc{\mA}{\mathbf{A}}
	\newc{\mB}{\mathbf{B}}
	\newc{\mC}{\mathbf{C}}
	\newc{\mD}{\mathbf{D}}
	\newc{\mE}{\mathbf{E}}
	\newc{\mF}{\mathbf{F}}
	\newc{\mJ}{\mathbf{J}}
	\newc{\mG}{\mathbf{G}}
	\newc{\mP}{\mathbf{P}}
	\newc{\mR}{\mathbf{R}}
	\newc{\mQ}{\mathbf{Q}}
	\newc{\mX}{\mathbf{X}}
	\newc{\muu}{\mathbf{u}}
	\newc{\mvv}{\mathbf{v}}
	
	\newc{\mllg}{\mathbb{\lambda}}
	\newc{\mLLg}{\mathbf{\LLg}}

	\newc{\lspn}[2]{\mbox{$\| #1\|_{\Lsp{#2}}$}}
	\newc{\Lpn}[2]{\mbox{$\| #1\|_{#2}$}}
	\newc{\Hn}[1]{\mbox{$\| #1\|_{H^1(\Og)}$}}

	\newc{\mynorm}[2]{\| #1\|_{#2}}

	\newcommand{\RR}{{\rm I\kern -1.6pt{\rm R}}}

	\newc{\itQQ}[2]{\dspl{\int_{#1}#2\,dz}}
	\newc{\mmitQQ}[2]{\dspl{\rule[1mm]{4mm}{.3mm}\hspace{-4.3mm}\int_{#1}~#2~dz}}
	\newc{\MmitQQ}[2]{\dspl{\rule[1mm]{4mm}{.3mm}\hspace{-4.3mm}\int_{#1}~#2~d\mu}}

	\newc{\MUmitQQ}[3]{\dspl{\rule[1mm]{4mm}{.3mm}\hspace{-4.3mm}\int_{#1}~#2~d#3}}
	\newc{\MUitQQ}[3]{\dspl{\int_{#1}~#2~d#3}}

	\newc{\mccP}{\mathbb{P}}
	\newc{\mccK}{\mathbb{K}}
	
	\newc{\DKTmU}{\mccK(U)}
	\newc{\DKTmUold}{(K_U(U)^{-1})^T}
	
	\newc{\myPi}{\mathbf{W}}
	\newc{\myIbar}{\bar{\ccI}_1}
	\newc{\myIhat}{\hat{\ccI}_1}
	\newc{\myIbreve}{\breve{\ccI}_0}
	
	\newc{\mmk}{\mathbf{k}}

	\newcommand{\ma}{\mathbf{a}}

	\newc{\mfu}{\mathbf{f_u}}
	\newc{\mh}{\mathbf{h}}
	\newc{\mb}{\mathbf{b}}
	\newc{\mf}{\mathbf{f}}

	\newcommand{\barrl}[2]{\barr{ll}\lefteqn{#1}\hspace{#2}&\\}

	\newc{\twomatrix}[1]{\left[\barr{cc}#1\earr\right]}
	\newc{\threematrix}[1]{\left[\barr{ccc}#1\earr\right]}

	\newc{\mN}{\mathbf{N}}
	\newc{\mI}{\mathbf{I}}
	\newc{\mH}{\mathbf{H}}

	\newc{\mk}{\mathbf{k}}
	\newc{\mr}{\mathbf{r}}

	\newc{\DIAGM}[2]{\left[\barr{ccc}#1&0\ldots&0\\
		\vdots&\ddots&\vdots\\0&\ldots0&#2\earr \right]}
	\newc{\DiagM}[2]{\mbox{diag}\left[#1
		\cdots #2 \right]}
	\newc{\vVEC}[2]{\left[\barr{c}#1\\
		\vdots\\#2\earr \right]}
	\newc{\hVEC}[2]{\left[#1
		\cdots #2 \right]}
	
	\newc{\mq}{\mathbf{q}}

	\newc{\msys}[1]{\left\{\barr{l}#1\earr
		\right.}
	\newc{\msysa}[1]{\left\{\barr{ll}#1\earr
		\right.}
	
	\newc{\bbM}{\mathbb{M}}
	\newc{\mat}[1]{\left[\barr{cc}#1\earr\right]}
	
	\newc{\me}{\mathbf{e}}

	\newc{\vecc}[2]{\left[\barr{cc}#1\\#2\earr\right]}
	\newc{\mL}{\mathbb{L}}
	
	\newc{\cO}{{\cal O}}
	
	\newc{\cM}{{\cal M}}
	
	\newc{\myega }{\eg_0(R)}
	\newc{\myeg}{\eg_1(\eg_*)}
	\newc{\myegp}{\hat{\eg}_1(\eg_*)}

\newc{\diagA}{\mathbb{A}_d}
\newc{\mBB}{\mathbb{B}}
\newc{\MLT}[1]{{\cal M}_{lt}(\Og,#1)}
\newc{\ALT}[1]{{\cal A}_{l}(\Og,#1)}
\newc{\mM}{\mathbb{M}}

\newc{\diag}[1]{\mbox{diag}(#1)}
\newc{\off}[1]{\mbox{offdiag}(#1)}
\newc{\mT}{\mathbb{T}}

\usepackage{amssymb}

\begin{document}

	\vspace*{-.8in}
	\begin{center} {\LARGE\em Global Existence to a Class of Triangular Block Matrix Cross Diffusion Systems and the Spectral Gap Condition}
		
	\end{center}

	\vspace{.1in}
	
	\begin{center}

		{\sc Dung Le}{\footnote {Department of Mathematics, University of
				Texas at San
				Antonio, One UTSA Circle, San Antonio, TX 78249. {\tt Email: Dung.Le@utsa.edu}
				}}

	\end{center}

	\begin{abstract} We study the  global existence of classical solutions  to cross diffusion systems of $m$ equations on  $N$-dimensional domains ($m,N\ge2$). The diffusion matrix is a triangular block matrix with coupled entries.
		We establish that the $W^{1,p}$ norm of solutions for some $p>N$ does not blow up in finite time so that the results in \cite{Am2} is applicable. We will also show that the spectral gap condition in \cite{dlebook,dlebook1} can be relaxed via a new result on BMO norms in \cite{dleBMO}. \end{abstract}

	\section{Introduction}\label{intro}
	
	This paper deals with the cross diffusion system
	\beqno{exsyszpara}\left\{\barr{ll}W_t=\Div(\ma(W)DW)+f(W)&\mbox{in $\Og\times (0,T)$,}\\ \mbox{homogeneous BC}&\mbox{on $\partial\Og\times(0,T)$,}\\W(x,0)=W_0(x)&\mbox{on $\Og$,  with $W_0\in W^{1,p_0}(\Og)$} \earr \right.\eeq
	on bounded domain $\Og$ of $\RR^N$ with homogeneous Dirichlet or Neumann boundary conditions. Here, $W$ is an unknown vector in $\RR^m$, $m\ge 2$. As usual, $\ma(W), f(W)$ are respectively $m\times m$ matrix and $m\times1$ vector. For some positive functions $\llg(W),\LLg(W)$ we assume the ellipticity condition
	\beqno{normalellcond}\LLg(W)|DW|^2\ge \myprod{\ma(W)DW,DW}\ge \llg(W)|DW|^2.\eeq

	The regularity of weak solutions to (regular or non-regular) \mref{exsyszpara} is a long standing problem in the theory of partial differential systems (e.g. see \cite{GiaS}). It is well known that  weak solutions are H\"older continuous around points where the BMO norms are small in small balls.	
	On the other hand, the global existence of strong/classical solutions has been considered in a pioneering paper of Amann \cite{Am2} where he showed that a solution would exist globally if $p_0>N$ and the $W^{1,p}(\Og)$ norm of $W$ for some $p>N$ does not blow up in finite time.
	
	Recently, in \cite{dlebook,dlebook1}, we introduce the strong/weak Gagliardo-Nirenberg  inequalities involving BMO norms (GNBMO for short) and assume that $\ma(W)$ satisfies  {\em the spectral gap condition} to prove that  the smallness of BMO norms in small balls
	\bdes \item[BMOsmall)] the BMO norm $\|W\|_{BMO(\mB_R)}$ on a ball $\mB_R\subset\Og$ is small if $R$ is small.
	\edes
	 also yields a control on the $W^{1,p}(\Og)$ norm of weak solutions and thus they are classical and exist globally.
	In \cite{dleBMO}  we  show that new weighted versions of the strong/weak GNBMO inequalities are also available and the crucial property BMOsmall) holds if $\ma(W)$ is regular elliptic.
	 This  allowed us to relax   the crucial partial regularity conditions (3.1) and (3.2) of \cite[Theorem 3.1]{GiaS} and establish that weak solutions of {\em regular} systems will be H\"older continuous just like those for regular scalar equations.

	We note that the spectral gap condition $$\nu_*=\sup_{W\in\RR^m}\frac{\llg(W)}{\LLg(W)}>1-2/N,$$ where $\llg(W),\LLg(W)>0$ are the smallest and largest eigenvalues of $\ma(W)$, has been used in \cite{dlebook1} is just technically sufficient for our theory. In fact, this condition was used only to guarantee that for some $c>0$ $$\myprod{\ma(W)DX,D(|X|^{2p-2}X)} \ge c|X|^{2p-2}|DX|^2, \quad X\in C^1(\Og,\RR^{m}) \mbox{ and } p>N/2$$ so that, together with BMOsmall) which holds for regular systems, we can derive an estimate for $\|DW\|_{L^{2p}(\Og)}$ and thus the global existence of strong solutions. We can see easily that the above inequality is true in special cases where the eigenvalue of $\ma(W)$ are far a part. For examples,  $\ma(W)$ can be triangular or a perturbation of a diagonal matrix with dominant diagonal entries.

	In this paper, we prove the global existence result for a nontrivial example where $\ma(W)$ can be {\em triangular block } matrix. Each diagonal sub-matrix satisfies the spectral gap condition but $\ma(W)$ does not.
	
	Consider the parabolic system on $\Og\times (0,T)$ ($\Og\subset\RR^N$ is a bounded domain, $N\ge2$, $T>0$) \beqno{trimat0}\left\{\barr{ll}W_t=\Div(\ma(W)DW)+f(W)&\mbox{in $\Og\times (0,T)$,}\\ \mbox{homogeneous BC}&\mbox{on $\partial\Og\times(0,T)$,}\\W(x,0)=W_0(x)&\mbox{on $\Og$,  with $W_0\in W^{1,p_0}(\Og)$.} \earr \right. \eeq Here, $W=(W_1,W_2)^T$,  $m=m_1+m_2$, $m_1,m_2\ge1$, $W_i\in\RR^{m_i}$ and  $$\ma(W)=\mat{A_1&A_3\\0&A_2},\quad f(W)=\mat{f_1\\f_2}.$$ $A_i$,'s are matrices and $f_i$'s are vectors.  $A_1,A_2$ are square matrices of sizes $m_1\times m_1$ and $m_2\times m_2$ respectively (so that $A_3$ is of size $m_1\times m_2$). $f_i$ is of size $m_i\times1$. The entries of $A_i,f_i$ are functions of $W$ (as we mentioned earlier, we can allow $\ma$ to depend on $x,t$ too as $W$ is assumed to be bounded), this assumption is crucial and makes our example of this paper nontrivial although $\ma(W)$ is not a completely full matrix.

	Our main result is
	
	\btheo{trimat-mainthm} Let $W$ be a {\em bounded} weak solution of \mref{trimat0} and $p_0>N$.	Assume that
	
	\bdes \item[i)] $\ma$ is smooth in $W$ and regular elliptic,
	\item[ii)] $A_1,A_2$ satisfy the spectral gap condition.
	\edes
	
	Then $W$ is classical and exists globally.
	
	\etheo
	
	Our paper is organized as follows. In \refsec{gensteps}, we sketch the results in \cite{dlebook1}  of proving the global existence.  We present the proof of \reftheo{trimat-mainthm} in \refsec{trimatsec} to conclude this paper. It is easy to see that our proof can cover the general case where $\ma(W)$ has $k\ge 2$ diagonal blocks and each of these satisfies spectral gap conditions.

	\section{Sketching of the main steps of proving the global existence in \cite{dlebook1}}\label{gensteps}\eqnoset
	In this section we recall the key ingredients of the proof of our existence theory and regularity properties of solutions to non-degenerate systems in  \cite{dlebook1}: the local weighted Gagliardo-Nirenberg inequalities involving BMO norm and its use in proving the global existence of \mref{exsyszpara}.

	\subsection{A new (weak) {\em weighted} Gagliardo-Nirenberg inequality} \index{(weak) Gagliardo-Nirenberg inequality} 
	Firstly, let us recall some notations.
	For any measurable subset $A$ of $\Og$  and any  locally integrable function $U:\Og\to\RR^m$ we denote by  $|A|$ the measure of $A$ and $U_A$ the average of $U$ over $A$. That is, $$U_A=\mitx{A}{U(x)} =\frac{1}{|A|}\iidx{A}{U(x)}.$$
	A function $f\in L^1(\Og)$ is said to be in $BMO(\Og)$  if for any ball  $Q\subset\Og$  $$[f]_{*}:=\sup_Q\mitx{Q}{|f-f_Q|}<\infty.$$ We then define $\|f\|_{BMO(\Og)}:=[f]_{*}+\|f\|_{L^1(\Og)}.$
	We assume that $H,Du$ are continuous and denote (here, $DH$ should be understood in its distribution sense). Fixing $R>0$, we let $\og$ be a smooth cutoff function in $B_R$ ($|D\og|\le C/R$).
	
	\newc{\mmyIbreve}{\mathbf{\myIbreve}}
	
	\beqno{Idefz} \mI_1:=\iidx{\Og}{|H|^{2p}|Du|^2\og^2},\;
	\mI_2:=\iidx{\Og}{|H|^{2p-2}|DH|^2\og^2},\;\mmyIbreve:=\iidx{\Og}{|H|^{2p}\og^2}.\eeq

	Especially, for a given $\eg_*>0$ we define
	\beqno{eg1def} \myeg=\left(\eg_*^{-N+2}\iidx{B_{\eg_*}}{|Du|^2}\right)^\frac12.\eeq

	We have the following inequality which will be referred to as the weak GNBMO weighted inequality for short throughout this paper. A more general (weighted) version of this inequality was reported in \cite[Theorem 2.1]{dleBMO}.
	\btheo{wGNBMOrecall}   Assume $\eg_*\le R$.
	Suppose that  $\myega, \mI_1,\mI_2,\mmyIbreve$ are finite.
	Then   there is a  constant $C$ such that \beqno{GNglobalestz}\mI_1\le \myega ^2 \mI_2+C(\frac{\myega ^2}{R^2}+\frac{\myeg^2}{\eg_*^2}+1)\mmyIbreve.\eeq
	\etheo

	\subsection{Regularity and global existence for the regular  system \mref{exsyszpara}:}
	However, assuming the {\em spectral-gap condition} and arguing as in \cite{dlebook1}, we see that under the assumption that the BMO norms of a weak (weak-strong) solution $W$ in small balls in $\Og$ are small then we can use the same induction argument in \cite[Lemma 4.4.1 or Theorem 4.4.6]{dlebook1} and the weak Gagliardo-Nirenberg BMO inequalities \reftheo{wGNBMOrecall} in finite steps to have a bound for $\|DW\|_{L^{2p}(B_R)}$ for some $p>N/2$. Of course, as $p>N/2$ then $W$ is H\"older continuous.
	
	We consider the systems
	\beqno{degsys} W_t=\Div(\ma(x,t,W)DW)+F(x,t,W), \eeq
	where there is $\llg(W),\LLg(W)>0$ such that $\llg(W)|\zeta|^2\le \myprod{\ma(x,t,W)\zeta,\zeta}\le \LLg(W)|\zeta|^2$. We drop $x,t$ in the sequel for simplicity. Since we will assume that $W$ is bounded, we see that the assumption that $D_x\ma(x,t,W)$ is bounded will be sufficient. Let $\nu_*=\sup_{W\in\RR^m}\frac{\llg(W)}{\LLg(W)}$.

	\btheo{regdegthm} Let  $W$ be a {\em bounded weak} solution to \mref{degsys}. {\em Suppose  the spectral-gap condition} $\nu_*>1-2/N$.  There is $p>N/2$ such that if $F_W\in L^{2p+2}(\Og)$.

	Then, for any $t_0>0$ there is a finite constant $C(p,t_0,\|W\|_{L^1(\Og)})$ such that $$\sup_{t\in(0,t_0)}\|DW\|_{L^{2p}(\Og_R)}\le C(p,t_0,\|F_W(W)\|_{L^{2p+2}(\Og)},\|W\|_{L^1(\Og)}).$$ Moreover, $W$ is H\"older continuous, $W\in C^\ag(\Og)$ for some $\ag>0$.
	
	Note that the {\em spectral-gap} is void if $m=1$ and we can have the result for all $\ag\in(0,1)$.

	\etheo
	
	The proof of \cite[Lemma 4.4.1]{dlebook1} to estimate $\|DW\|_{L^{2p}_{loc}(\Og)}$ is as follows. Applying the difference operator  $\dg_h$ in $x$ to  the equation of $W$, we see that $W$ weakly solves (dropping the terms involving $\mB$ in the proof of \cite[Lemma 4.4.1]{dlebook1})
	$$(\dg_hW)_t=\Div(\ma D(\dg_h W)+\ma_W \myprod{\dg_h W,DW})+F_W  \dg_hW.$$
	
	Let $p\ge1$ and   $\fg,\eta$ be positive $C^1$ cutoff functions for the concentric balls $B_s, B_t$ and the time interval $I$. For any $0<s<t<2R_0$ we test this system with $|\dg_hW|^{2p-2}\dg_hW\fg^2\eta$ and use Young's inequality for the term $|\ma_W||\dg_hW|^{2p-1}|DW||D(\dg_hW)|$ and the spectral gap condition (with $X=\dg_h W$). We get    for some constant $c_0>0$, $\Psi=|F_W|$,  and $Q=\Og_t\times I$
	$$\barrl{\sup_{t\in(t_0,T)}\iidx{\Og}{|\dg_hW|^{2p}\fg^2}+c_0\itQ{Q}{ \llg|\dg_hW|^{2p-2}|D(\dg_hW)|^2\fg^2}\le  }{1cm}& C\iidx{\Og\times\{t_0\}}{|\dg_hW|^{2p}\fg^2}+C\itQ{Q}{\frac{|\ma_W|^2}{\llg}|\dg_hW|^{2p}|DW|^2\fg^2}+\\&C\itQ{Q}{[\Psi|\dg_hW|^{2p}\fg^2+|\ma_W|(|\dg_hW|^{2p}|DW|)\fg|D\fg|]}.\earr$$
	
	We make use of \reftheo{wGNBMOrecall}, with $H=Du=\dg_h W$,   to estimate the second term on the right hand side and absorb it into the left hand side ($W$ is bounded and so is $\frac{|\ma_W|^2}{\llg}$). Here, we make use of a very important result in \cite{dleBMO} which established 
BMOsmall) in the Introduction.
	
	This result holds because we were able to establish that: 
	{\em Higher integrability in spatial derivatives $\Rightarrow$ BMOsmall) in space and time.} We  established new theoretic functional results in \cite{dleBMO}.
	One of the main themes of \cite[section 5.3]{dleBMO} is the applications of these new theoretic functional results in combining with the higher integrabilities of spatial derivatives in \cite{GiaS} : If \mref{degsys} is {\em regular elliptic}, then there are $p_*>2$ and a constant $C$ such that  $DW\in L^{p_*}_{loc}(\Og\times (0,T))$ and for any $R>0$ and $B_R\subset\Og$, via translations,
	$$\left(\frac{1}{R^{d+2}}\dspl{\int_0^{R^2}}\iidx{B_R}{|DW|^{p_*}}ds\right)^\frac{1}{p_*}\le C \left(\frac{1}{R^{d+2}}\dspl{\int_0^{R^2}}\iidx{B_R}{|DW|^{2}}ds\right)^\frac{1}{2}$$
	to establish that a weak solution $W$ to  parabolic cross diffusion systems on a domain $\Og\subset\RR^d$ (respectively $\Og\times(0,T)$, $T>0$) with $d\ge2$ has the BMOsmall) property.

	We then obtain the estimate for $\sup_I\|DW\|_{L^{2p}(\Og_R)}$ for some $p>N/2$. We will modify these ideas to obtain the proof of our main result in the next section.

\section{The proof of our main theorem} \label{trimatsec}\eqnoset We now present the proof of our main result and note that the spectral gap condition $\nu_*>1-2/N$  is just technically sufficient for our theory. In fact, this condition was used only to guarantee that for some $c>0$ $$\myprod{\ma(W)DX,D(|X|^{2p-2}X)} \ge c|X|^{2p-2}|DX|^2, \quad X\in C^1(\Og,\RR^{m}) \mbox{ and } p>N/2$$ so that, together with BMOsmall) which holds for regular systems, we can derive an estimate for $\|DW\|_{L^{2p}(\Og)}$ and thus the global existence of strong solutions. We can see easily that the above inequality is true in special cases where the eigenvalue of $\ma(W)$ are far a part. For examples,  $\ma(W)$ can be triangular or a perturbation of a diagonal matrix with dominant diagonal entries.

In this paper, we prove the global existence result for a nontrivial example where $\ma(W)$ can be {\em triangular block } matrix. Each diagonal sub-matrix satisfies the spectral gap condition but $\ma(W)$ does not.

For the convenience of readers, we repeat the statement of our main result (\reftheo{trimat-mainthm}).
Consider the parabolic system on $\Og\times (0,T)$ ($\Og\subset\RR^N$ is a bounded domain, $N\ge2$, $T>0$) \beqno{trimat}\left\{\barr{ll}W_t=\Div(\ma(W)DW)+f(W)&\mbox{in $\Og\times (0,T)$,}\\ \mbox{homogeneous BC}&\mbox{on $\partial\Og\times(0,T)$,}\\W(x,0)=W_0(x)&\mbox{on $\Og$,  with $W_0\in W^{1,p_0}(\Og)$.} \earr \right. \eeq Here, $W=(W_1,W_2)^T$,  $m=m_1+m_2$, $m_1,m_2\ge1$, $W_i\in\RR^{m_i}$ and  $$\ma(W)=\mat{A_1&A_3\\0&A_2},\quad f(W)=\mat{f_1\\f_2}.$$ 
	
\btheo{trimatcoro} Let $W$ be a {\em bounded} weak solution of \mref{trimat} and $p_0>N$.	Assume that
	
\bdes \item[i)] $\ma$ is smooth in $W$ and regular elliptic,
\item[ii)] $A_1,A_2$ satisfy the spectral gap condition.
\edes

Then $W$ is classical and exists globally.

\etheo

Note that the assumption $W$ is bounded can be realized easily by \cite[Theorem 5.10]{dleBMO} under appropriate conditions for $f(W)$. 

\bproof
The assumption i) implies that  a higher integrability of $DW$ is available for \mref{trimat} so that  $W$ is H\"older continuous and  $W$ satisfies BMOsmall), see \cite{dleBMO}. That is \beqno{BMOtrimat} \|W_1\|_{BMO(\Og_R)} , \|W_2\|_{BMO(\Og_R)}\mbox{ are small if $R$ is sufficiently small.}\eeq
	
Explicitly,	\mref{trimat} implies two subsystems \beqno{trimat1}(W_1)_t=\Div(A_1(W_1,W_2)DW_1+A_3(W_1,W_2)DW_2)+f_1(W_1,W_2),\eeq
	\beqno{trimat2}(W_2)_t=\Div(A_2(W_1,W_2)DW_2)+f_2(W_1,W_2).\eeq
	
The plan of our proof is simple, we study the regularity of $W_2$ (solution to the system\mref{trimat2} which is still coupled in $A_2,f_2$ with \mref{trimat1}). This regularity is then used in \mref{trimat1}. This is nontrivial as $A_2$ depends on $W_1$.

By ii), we can find $c>0$ and $p>N/2$ such that $$\myprod{A_i(W)DX,D(|X|^{2p-2}X)} \ge c|X|^{2p-2}|DX|^2, \quad X\in C^1(\Og,\RR^{m}),\; i=1,2.$$

We follow the proof of \cite[Lemma 4.4.1]{dlebook1} to estimate $\|DW_2\|_{L^{2p}_{loc}(\Og)}$ as in we did in \reftheo{regdegthm} for $\|DW\|_{L^{2p}_{loc}(\Og)}$. We differentiate (more precisely, take the difference) \mref{trimat2} to have (with a slight abuse of notation)
\beqno{trimat2a}(DW_2)_t=\Div(A_2D^2W_2+(A_2)_{W_2}DW_2DW_2+(A_2)_{W_1}DW_1DW_2)+(f_2)_{W_2}DW_2+(f_2)_{W_1}DW_1.\eeq 

Here we have the extra terms  $(A_2)_{W_1}DW_1DW_2$ and $(f_2)_{W_1}DW_1$. Testing the above with $|DW_2|^{2p-2}DW_2\fg^2$ and using ii) and  the facts that $(A_2)_{W_1}, (f_2)_{W_1}$ are bounded because $W$ is H\"older continuous, we obtain for $Q=\Og\times I$, where $I$ is any time interval and $\fg(x)$ is a cut-off function, the following estimate for some constant $C,c>0$
$$\barrl{\sup_I\iidx{\Og}{|DW_2|^{2p}\fg^2}+c\itQ{Q}{|DW_2|^{2p-2}|D^2W_2|^2\fg^2}\le C(1/R)\itQ{Q}{|DW_2|^{2p}\fg^2}+}{6cm}&C\itQ{Q}{(|DW_2|^{2p}|DW_1|^2+|DW_2|^{2p-1}|DW_1|)\fg^2}.\earr$$
Here, the extra terms  $(A_2)_{W_1}DW_1DW_2$ and $(f_2)_{W_1}DW_1$ give rise to the last integral. Using $$|DW_2|^{2p-1}|DW_1|=|DW_2|^{2p-2}|DW_1||DW_2|\le C(|DW_2|^{2p-2}|DW_1|^2+|DW_2|^{2p},$$  we derive some other constants $C,c>0$ such that
\beqno{trimatkey}\barrl{\sup_I\iidx{\Og}{|DW_2|^{2p}\fg^2}+c\itQ{Q}{|DW_2|^{2p-2}|D^2W_2|^2\fg^2}\le}{4cm}& C(1/R)\itQ{Q}{|DW_2|^{2p}\fg^2}+C\itQ{Q}{|DW_2|^{2p}|DW_1|^2\fg^2}.\earr\eeq

To estimate the last term, we make use of the weak GNBMO inequality \cite[Theorem 2.4.2]{dlebook1}
\beqno{weakGNBMOtrimat}\iidx{\Og}{|H|^{2p}|Du|^2\fg^2}\le \eg_0^2(R)\iidx{\Og}{|H|^{2p-2}|DH|^2\fg^2}+R^{-2}\eg_1^2(R)\iidx{\Og}{|H|^{2p}\fg^2}\eeq
	where $R>0$, $\fg$ is a cut-off function for $B_R,B_{2R}$ and $$\eg_0(R)=\|u\|_{BMO(R)},\; \eg_1(R)=\left(R^{-N+2}\iidx{B_{2R}}{|Du|^2}\right)^\frac{1}{2}.$$
	
We  let $H=DW_2, u=W_1$ and fix $R$ small such that $\eg_0$ is sufficiently small. Then, the last term in \mref{trimatkey} can be absorbed to the second one on the left hand side. We then see that $\|DW_2\|_{L^{2p}(\Og)}$ can not blow up in finite time. As $W$ is H\"older continuous then so is $A_2$, \cite[Theorem 3.2]{GiaS} gives that $DW_2$ is H\"older locally continuous (in $x$) with some exponent $\sg>0$. This then shows that \beqno{A3DW2}B(x,t):=A_3(W)DW_2 \mbox{ satisfies BMOsmall).}\eeq

Now, we turn to\mref{trimat1} which can be written as
$$(W_1)_t=\Div(A_1DW_1 + B)+f_1(W).$$

Because $A_1$ satisfies a spectral gap condition, we can follow the proof of \cite[Lemma 4.4.1]{dlebook1} again to have an estimate for $\|DW_1\|_{L^{2p}(\Og)}$. That is, we  We differentiate (more precisely, take the difference) the above and test the result with $|DW_1|^{2p-2}DW_1\fg^2$ to have for some constants $C,c>0$
$$\barrl{\sup_I\iidx{\Og}{|DW_1|^{2p}\fg^2}+c\itQ{Q}{|DW_1|^{2p-2}|D^2W_1|^2\fg^2}\le}{5cm}& C(1/R)\itQ{Q}{|DW_1|^{2p}\fg^2}+C\itQ{Q}{|DW_1|^{2p-2}|DB|^2\fg^2}.\earr$$

Again, we use \mref{weakGNBMOtrimat} ($p$ is now $p-1$) and  let $H=DW_1, u=B$ and fix $R$ small such that $\eg_0$ is sufficiently small, thanks to \mref{A3DW2}. Then, the last term in the above estimate can be absorbed to the second one on the left hand side. We then see that $\|DW_1\|_{L^{2p}(\Og)}$ can not blow up in finite time.

Thus, we show that $\|DW_1\|_{L^{2p}(\Og)}$, $\|DW_2\|_{L^{2p}(\Og)}$ do not blow up in finite time and complete the proof by using the result of \cite{Am2}.
\eproof

\brem{BMOtrimatrem} We see the crucial role of the result in \cite{dleBMO} showing that higher integrability of gradients implies BMOsmall). The system \mref{trimat} is not strongly coupled because $\ma$ is not a full matrix. However, as the entries of $A_2$ depend on $W_1$ so that the system is not totally weakly coupled and thus proving that $W_2$ is not a simple matter as the assumed boundedness of $W$ is not sufficient. Thanks to \cite{dleBMO} we have that $W$ satisfies BMOsmall)  and this is essential in estimating the last term of \mref{trimatkey} in order to prove the global existence of $W_2$ and then $W_1$.

\erem

\brem{advec} It is easy to see that the proof is similar if we  add an advection term $\mb(W)DW$  to the system \mref{trimat} and consider
$$W_t=\Div(\ma(W)DW)+\mb(W)DW+f(W) \mbox{ with }\mb(W)=\mat{B_1&B_3\\0&B_2}.$$ 
 The entries of $B_i$'s must be  vectors.

\erem

	\bibliographystyle{plain}

\end{document}